\documentclass[a4paper,12pt]{amsart}

\usepackage{amssymb, manfnt}
\usepackage{hyperref}
\usepackage[normalem]{ulem}

\def\altdb{\vadjust{\vbox to 0pt{\vss\hbox{\kern \hsize
\quad{\dbend}}\kern\baselineskip\kern-10pt}}}

\newcommand\field[1]{\mathbb{#1}}
\newcommand\CC{\field{C}}
\newcommand\NN{\field{N}}

\newcommand\ZZ{\field{Z}}

\newcommand\Gg{G}

\newcommand\Oo{\mathcal O}

\newcommand\Uu{\mathcal U}

\renewcommand\ker{\operatorname{ker}}

\newcommand\id{\operatorname{id}}

\newcommand\lsp{\operatorname{span}}

\newcommand\supp{\operatorname{supp}}

\newcommand{\ag}{A(\Gg)}
\newcommand{\agL}{A(\Gg_\Lambda)}
\newcommand{\Gtwo}{\Gg^{(2)}}

\newcommand{\cs}{\ensuremath{C^{*}}}

\theoremstyle{plain}
\newtheorem{theorem}{Theorem}[section]
\newtheorem*{theorem*}{Theorem}
\newtheorem*{prop*}{Proposition}
\newtheorem{cor}[theorem]{Corollary}
\newtheorem{lemma}[theorem]{Lemma}
\newtheorem{prop}[theorem]{Proposition}

\theoremstyle{remark}
\newtheorem{rmk}[theorem]{Remark}

\newtheorem{example}[theorem]{Example}

\theoremstyle{definition}
\newtheorem{dfn}[theorem]{Definition}

\newtheorem{notation}[theorem]{Notation}

\newcommand{\Gp}{\Gamma}
\newcommand{\go}{\Gg^{(0)}}
\newcommand{\Bco}[2]{B^{\operatorname{co}}_{#1}(#2)}
\newcommand{\BcoG}{\Bco{*}{\Gg}}
\newcommand{\BncoG}{\Bco{n}{\Gg}}

\usepackage{color}

\numberwithin{equation}{section}

\title{A groupoid generalisation of Leavitt path algebras}

\author{Lisa Orloff Clark}
\email{lclark@maths.otago.ac.nz}
\address{Lisa Orloff Clark\\ Department of Mathematics and Statistics\\ University of Otago\\ PO Box 56\\ Dunedin 9054\\ NEW ZEALAND}

\author{Cynthia Farthing}
\email{cynthia-farthing@uiowa.edu}
\address{Cynthia Farthing\\ Department of Mathematics\\ 14 MacLean Hall\\ Iowa City\\ Iowa 52242-1419\\ USA}

\author{Aidan Sims}
\email{asims@uow.edu.au}
\address{Aidan Sims\\ School of Mathematics and Applied Statistics\\ University of Wollongong\\ Wollongong\\ NSW 2518 \\ AUSTRALIA}

\author{Mark Tomforde}
\email{tomforde@math.uh.edu}
\address{Mark Tomforde\\ Department of Mathematics\\ University of Houston\\ Houston\\ TX 77204-3008\\ USA}

\date{\today}
\keywords{Topological groupoids; Leavitt algebra; groupoid algebra; graded algebra}
\subjclass[2010]{16S99 (Primary); 16S10, 22A22 (Secondary)}

\begin{document}

\begin{abstract}
Let $G$ be a locally compact, Hausdorff groupoid in which $s$ is a local homeomorphism and $\go$ is
totally disconnected.  Assume there is a continuous cocycle $c$ from $G$ into a discrete group
$\Gamma$. We show that the collection $A(G)$ of locally-constant, compactly supported
functions on $G$ is a dense $*$-subalgebra of $C_c(G)$ and that it is universal for algebraic
representations of the collection of compact open bisections of $G$. We also show that if $G$ is
the groupoid associated to a row-finite graph or $k$-graph with no sources, then $A(G)$ is
isomorphic to the associated Leavitt path algebra or Kumjian-Pask algebra. We prove versions of the
Cuntz-Krieger and graded uniqueness theorems for $A(G)$.
\end{abstract}

\maketitle

\section{Introduction}\label{sec:intro}

A ring $R$ is said to have invariant basis number if any two bases (i.e., $R$-linearly independent
spanning sets) of a free left $R$-module have the same number of elements.  Many familiar rings
(e.g., fields, commutative rings, left-Noetherian rings) have invariant basis number, but there are
many examples of noncommutative rings that do not.  A ring $R$ without invariant basis number is
said to have module type $(m,n)$ if $m < n$ are {natural numbers} chosen minimally with $R^m \cong
R^n$ as left $R$-modules.  In the 1940's, Leavitt constructed algebras $L_{m,n}$ with module type
$(m,n)$ for all pairs of natural numbers with $m < n$ \cite{Leav2, Leav3}.  The $L_{m,n}$ are now
known as the \emph{Leavitt algebras}, and when $m = 1$, the Leavitt algebra $L_{1,n}$ is the unique
nontrivial unital complex algebra generated by elements $x_1 \dots x_n$ and $y_1, \dots, y_n$ such
that $\sum^n_{i=1} x_i y_i = 1$ and $y_i x_j = \delta_{i,j} 1$ for all $i,j \le n$. In the 1970's,
independent of Leavitt's work and motivated by the search for $\cs$-algebraic analogues of Type~III
factors, Cuntz defined a class of $\cs$-algebras $\Oo_n$, one for each integer $n \ge 2$, which are
generated by elements $s_1, \dots, s_n$ satisfying $\sum^n_{i=1} s_i s^*_i = 1$ and $s^*_i s_i = 1$
for all $i$ (it follows that $s^*_i s_j = \delta_{i,j} 1$ for all $i,j \le n$). A consequence of
the uniqueness of $L_{1,n}$ is that it is isomorphic to the dense $*$-subalgebra of $\Oo_n$
generated by $s_1, \dots, s_n$ via an isomorphism that carries each $x_i$ to $s_i$ and each $y_i$
to $s^*_i$.

Shortly after Cuntz's work, Cuntz and Krieger generalised Cuntz's results to describe a class of
$C^*$-algebras $\Oo_A$ associated to binary-valued matrices $A$ \cite{CuntzKrieger:IM80}. At about
the same time, Enomoto and Watatani provided a very elegant description of these Cuntz-Krieger
algebras in terms of the directed graphs encoded by the matrices. Nearly twenty years later,
Kumjian, Pask, Raeburn, and Renault developed the class of $C^*$-algebras now known as graph
$C^*$-algebras \cite{KumjianPaskEtAl:JFA97}, as a far-reaching generalisation of the Cuntz-Krieger
algebras patterned on Enomoto and Watatani's approach. Each graph $C^*$-algebra is described in
terms of generators associated to the vertices and edges in the graph subject to relations encoded
by connectivity in the graph. The Cuntz algebra $\Oo_n$ corresponds to the graph with one vertex
and $n$ edges. A remarkable assortment of important $C^*$-algebraic properties of a graph
$C^*$-algebra can be characterised in terms of the structure of the graph (see \cite{Raeburn:CMBS}
for a good overview). Shortly afterwards, Kumjian and Pask introduced a sort of higher-dimensional
graph \cite{KumjianPask:NYJM00}, now known as a $k$-graph, and an associated class of
$C^*$-algebras, as a flexible visual model for the higher-rank Cuntz-Krieger algebras discovered by
Robertson and Steger \cite{RobertsonSteger:JRAM99}. When $k = 1$, a $k$-graph is essentially a
directed graph, and Kumjian and Pask's $C^*$-algebras coincide with the graph $C^*$-algebras of
\cite{KumjianPaskEtAl:JFA97}.

In the early 2000's, the algebraic community became interested in the similarity between the
constructions of Leavitt and Cuntz and the potential for the graph $C^*$-algebra template to
provide a broad class of interesting new algebras. Following the lead of
\cite{KumjianPaskEtAl:JFA97}, Abrams and Aranda Pino associated \emph{Leavitt path algebras} to a
broad class of directed graphs. The Leavitt path algebra of a directed graph is the universal
algebra whose presentation in terms of generators and relations is essentially the same as that of
the graph $C^*$-algebra. Moreover, the graded uniqueness theorem for Leavitt path
algebras implies that the $C^*$-algebra of a directed graph is a norm completion of its Leavitt
path algebra \cite{Malaga}, \cite{Tom10}. Further generalising Leavitt path algebras, Aranda Pino,
(J.)~Clark, an Huef, and Raeburn recently constructed a class of algebras associated to $k$-graphs,
which they call Kumjian-Pask algebras \cite{ArandaClarkEtAl:xx11}.

A very powerful framework for constructing $C^*$-algebras is the notion of a groupoid
$C^*$-algebra. Renault's structure theory for groupoid $C^*$-algebras \cite{Ren} is exploited in
\cite{KumjianPaskEtAl:JFA97} where structural properties of the graph $\cs$-algebra are deduced by
showing that the graph $\cs$-algebra is isomorphic to a groupoid $C^*$-algebra and then tapping
into Renault's results \cite{Ren}.
The same approach was taken in \cite{KumjianPask:NYJM00} to establish important structural
properties of $k$-graph $C^*$-algebras: the $C^*$-algebra of a $k$-graph is defined in terms of
generators and relations, but its structure is analysed by identifying it with a groupoid
$C^*$-algebra.

In this paper, from a sufficiently well-behaved groupoid $\Gg$, we construct a
complex algebra $\ag$ with
the following properties:
\begin{enumerate}
 \item $\ag$ has a natural description as a universal algebra
     (Theorem~\ref{thm:UniversalProp});
    \item $\ag$ is isomorphic to a dense subalgebra of the groupoid $\cs$-algebra $C^*(\Gg$)
        (Proposition~\ref{prop:AG dense in CcG}); and
    \item {given a $k$-graph $\Lambda$, if  $\Gg = \Gg_\Lambda$ is the groupoid corresponding
        to $\Lambda$ as in \cite{KumjianPask:NYJM00} (Proposition~\ref{prop:agL isomorphic to
        KP algebra}), then $\ag$ is isomorphic to the Kumjian-Pask algebra
        $\operatorname{KP}_\CC(\Lambda)$.} In particular, if $E$ is a directed graph and $\Gg =
\Gg_E$ is the graph groupoid associated to $E$, then $A(\Gg)$ is isomorphic to the Leavitt path
algebra $L_\mathbb{C}(E)$.
\end{enumerate}

In \cite{Steinberg2010}, Steinberg  defines a groupoid
algebra $KG$ for an arbitrary commutative ring $K$ with unit
 and shows that $KG$ is isomorphic to
 an associated  inverse semigroup algebra.
We show that the algebra $\ag$ is identical to $KG$ for $K = \CC$
(the complex numbers).\footnote{We would like to thank Steinberg who brought this to our attention
after reading an earlier version of this paper.}
Our approach is significantly different from that of \cite{Steinberg2010} and
our universal-property and uniqueness theorems (see below) provide tools for studying $KG$
and the inverse semigroup algebras associated to them in the case where  $K= \CC$.
Conversely,  it would be interesting to investigate
versions of our theorems for general $K$.

The Cuntz-Krieger uniqueness theorem and gauge-invariant uniqueness theorem are important tools in
the study of graph $C^*$-algebras.  Versions of these theorems have been established for many
generalisations of Cuntz-Krieger algebras \cite{BPRS,FMY, KumjianPask:NYJM00, KPR1998,
KumjianPaskEtAl:JFA97, RSY1, RSY2}. For Leavitt path algebras, the graded uniqueness theorem is the
analogue of the gauge-invariant uniqueness theorem.  The first version of this graded
uniqueness theorem was a corollary to Ara, Moreno, and Pardo's characterisation
\cite[Theorem~4.3]{AMP} of the graded ideals in a Leavitt path algebra. It was first stated
explicitly by Raeburn who proved both the graded uniqueness theorem and Cuntz-Krieger uniqueness
theorem for Leavitt path algebras of row-finite graphs with no sinks and over fields equipped with
a positive definite $*$-operation \cite[Theorem 1.3.2 and Theorem 1.3.4]{Malaga}. Tomforde extended
these results to Leavitt path algebras of arbitrary graphs over arbitrary fields in
\cite[Theorem~4.8 and Theorem~6.8]{Tom10}, and later proved the two uniqueness theorems for Leavitt
path algebras of arbitrary graphs over a ring \cite[Theorem 5.3 and Theorem 6.5]{Tom12}.  Aranda
Pino, (J.) Clark, an Huef, and Raeburn subsequently proved versions of these theorems for
Kumjian-Pask algebras \cite{ArandaClarkEtAl:xx11}. In Section~\ref{sec:uniqueness} we prove
versions of the Cuntz-Krieger uniqueness theorem (Theorem \ref{thm:CKUT}) and the graded uniqueness
theorem (Theorem \ref{thm:kernel of graded homomorphism}) for $A(G)$. We also give an example of a
groupoid satisfying our hypothesis that is not necessarily the groupoid of a $k$-graph.

Our aim in defining and initiating the analysis of $A(G)$ is twofold: (1) to provide a broad
framework for future generalisations of Leavitt path algebras from other combinatorial structures;
and (2) to make available the powerful toolkit of groupoid analysis to study these algebras. In
addition, we hope this will provide a new and useful perspective on the interplay between algebra
and analysis at the interface between  Leavitt path algebras and graph $C^*$-algebras.

\section{Preliminaries}\label{sec:prelims}

A groupoid is a small category with inverses. We  write $\Gtwo \subseteq \Gg \times \Gg$ for the
set of composable pairs in $\Gg$; we  write $\go$ for the unit space of $\Gg$, and we  denote by
$r$ and $s$ the range and source maps $r,s : \Gg \to \go$. So $(\alpha,\beta) \in \Gtwo$ if
$s(\alpha) = r(\beta)$. For $U,V \subseteq G$, we define
\begin{equation}\label{eq:setproduct}
UV:= \{ \alpha \beta : \alpha \in U, \ \beta \in V,
        \text{ and } r(\beta) = s(\alpha) \}.
\end{equation}
A topological groupoid is a groupoid endowed with a topology under which $r$ and $s$ are
continuous, the inverse map is continuous, and such that composition is continuous with respect to
the relative topology on $\Gtwo$ inherited from $\Gg \times \Gg$.

Recall that if $\Gg$ is a groupoid, then an \emph{open bisection} of $\Gg$ is an open subset $U
\subseteq \Gg$ such that $r|_U$ and $s|_U$ are homeomorphisms.
We will work exclusively with locally compact,
Hausdorff groupoids in which the source map $s : \Gg \to \go$ is a local homeomorphism.

\begin{rmk}\label{rmk:Renault etale}
In modern nomenclature, a groupoid in which the source map is a local homeomorphism is said to be
\'etale. The range map is then a local homeomorphism as well since inversion is a continuous,
self-inverse bijection that interchanges $s$ and $r$. Whenever $s : \Gg \to \go$ is a
local homeomorphism, it follows that $\go$ is open in $\Gg$, which in modern terminology is to say
that $\Gg$ is \emph{$r$-discrete}. The converse does not hold without additional hypotheses: if
$\go$ is open in $\Gg$, then $s$  is a local homeomorphism if and only if $\Gg$ admits a left Haar
system, in which case the Haar system consists of counting measures (see the opening paragraph of
Section~3 in \cite{QuiggSieben:JAMSSA99}). In the first papers on the subject, the term \'etale
referred to a groupoid in which $\go$ is open, and what is now known as an \'etale groupoid was
described as an \'etale groupoid with a Haar system.  The terminology is now a little ambiguous.
We will circumvent the issue by simply writing out what we are assuming about our groupoids. Since,
in particular, in \cite{Ren} the phrase ``an \'etale groupoid with a Haar system" is equivalent to
``a groupoid in which the source map is a local homeomorphism," we can and will apply results from
\cite{Ren}.
\end{rmk}

\begin{lemma}\label{lemma:clopen bisection basis}
Let $\Gg$ be a topological groupoid such that $s : \Gg \to \go$ is a local homeomorphism. Suppose
that $\go$ has a basis consisting of clopen sets; that is, $\go$ is totally disconnected. Then
the topology on $\Gg$ has a basis of clopen bisections. Moreover, if $\Gg$ is locally
compact and Hausdorff, then $\Gg$ has a basis of compact open bisections.
\end{lemma}
\begin{proof}
By Remark~\ref{rmk:Renault etale}, we may apply Proposition~2.8 of \cite{Ren} to see that $\Gg$ has
a basis of open bisections. For each $\gamma \in \Gg$, let $U$ be an open bisection containing
$\gamma$. Since $r$ is an open map there exists a basic clopen neighbourhood $X$ of $r(\gamma)$
such that $X \subseteq r(U)$. Then $XU = \{h \in U : r(h) \in X\} = U \cap r^{-1} (X)$ is
homeomorphic to $X$ by choice of $U$ and in particular is a clopen bisection containing $\gamma$.
If $\Gg$ is also locally compact, then $U$ may be chosen to be precompact. Hence the clopen subset
$XU$ is a compact open bisection.
\end{proof}

\begin{notation}
For the remainder of this paper, $\Gp$ will denote a discrete group, $\Gg$ will denote a locally
compact, Hausdorff groupoid with totally disconnected unit space in which $s : \Gg \to \go$ is a
local homeomorphism, and $c$ will denote a continuous cocycle from $\Gg$ to $\Gp$ (that is, $c$
carries composition in $\Gg$ to the group operation in $\Gp$).
\end{notation}

By Lemma~\ref{lemma:clopen bisection basis}, with $\Gp$, $\Gg$ and $c$ as above, $\Gg$ has a basis
of compact open bisections. Since $\Gg$ is Hausdorff, compact subsets of $\Gg$ are closed. We will
use this fact frequently and without further comment.

\begin{rmk}
These hypotheses might sound very restrictive, but, for instance, every $k$-graph groupoid satisfies
them (see, for example, \cite{FMY}).
\end{rmk}

\begin{rmk}\label{rmk:disjointify}
Let $U$ be a compact open subset of a topological space $X$. Let $F$ be a finite cover of $U$ by
compact open subsets of $U$. For each nonempty $H \subseteq F$, let $V_H := \big(\bigcap H\big)
\setminus \big(\bigcup (F \setminus H)\big)$. Since each $V \in F$ is compact and open, so is each
$V_H$.  In particular, since $F$ is finite, so is $K := \{H \subseteq F : H \not= \emptyset, V_H
\not= \emptyset\}$, and
\[
U = \bigsqcup K
\]
is an expression for $U$ as a finite disjoint union of  nonempty compact open sets such that
for each $W \in K$ we have $W \subseteq V$ for at least one $V \in F$, and such that
whenever $W \in K$ and $V \in F$ satisfy $W \not\subseteq V$, we have $W \cap V =
\emptyset$. We refer to this as the \emph{disjointification} of the cover $F$ of $U$.
\end{rmk}

\section{Construction of the Leavitt groupoid algebra}\label{sec:construction}

\begin{dfn}
Let $X$ be a topological space.  A function $f : X \to Y$ is
\emph{locally constant} if for every $x \in X$ there exists a
neighbourhood $U$ of $x$ such that $f|_U$ is constant.
\end{dfn}

\begin{lemma}\label{lem:clopen support}
If $f : X \to \CC$ is continuous and locally constant then $\supp(f) := \{x \in X : f(x) \not= 0\}$ is clopen.
\end{lemma}
\begin{proof}
Since $f$ is locally constant, for each $x \in X$ such that $f(x) = 0$, there is an open
neighbourhood $U^0_{x}$ of $x$ such that $f(y) = 0$ for all $y \in U^0_{x}$. Hence $\{x \in X : f(x)
= 0\} = \bigcup_{f(x) = 0} U^0_x$ is open. It is also closed as it is the preimage of the closed
set $\{0\}$. Hence it and its complement $\supp(f)$ are both clopen.
\end{proof}

\begin{dfn}
Let $\Gp$ be a discrete group, $\Gg$ a locally compact, Hausdorff groupoid with totally
disconnected unit space such that $s : \Gg \to \go$ is a local homeomorphism, and $c : \Gg \to \Gp$
a continuous cocycle. For each $n \in \Gp$, let $\Gg_n := c^{-1}(n) \subseteq \Gg$, and let $A_n$
denote the complex vector space
\[
    A_n := \{f\in C_c(\Gg) \ | \ \supp(f) \subseteq \Gg_n \text{ and $f$ is locally constant} \}
\]
with pointwise addition and scalar multiplication. We define $\ag := \lsp\{A_n : n \in \Gp\}
\subseteq C_c(\Gg)$. The $A_n$ are pairwise linearly independent in
$\ag$.
\end{dfn}

By definition, $\ag$ consists of compactly-supported functions $f : \Gg \to \CC$ such that each
$f_n := f|_{\Gg_n}$ belongs to $A_n$, and $f_n = 0$ for all but finitely many $n$.
Lemma~\ref{lem:clopen support} implies that $f \in C_c(\Gg)$ is locally constant if and only if
$f|_{\supp(f)}$ is locally constant.

The notation $\ag$ does not suggest any dependence on the cocycle $c$. The following lemma
justifies this. Since we may endow any groupoid $\Gg$ with the trivial cocycle into the trivial
group $\{e\}$, it follows that $\ag$ is defined for any groupoid $\Gg$ with totally disconnected,
locally compact unit space such that $s$ is a local homeomorphism.   This lemma then
verifies that $A(G)$ is precisely the algebra $\CC G$ of
\cite[Definition~4.1]{Steinberg2010}.

\begin{lemma}\label{lem:span}
Let $\Gp$ be a discrete group, $\Gg$ a locally compact, Hausdorff groupoid with totally disconnected
unit space such that $s : \Gg \to \go$ is a local homeomorphism, and $c : \Gg \to \Gp$ a continuous
cocycle. Then $\ag = \{f \in C_c(\Gg) : f\text{ is locally constant}\}$. If $\Uu$ is the basis of
all compact open subsets of $\Gg$, we have $\ag = \lsp\{1_U : U \in \Uu\}$.
\end{lemma}
\begin{proof}
Let $f \in \ag$. Write $f = \sum_{k \in F} f_k$ where $F \subseteq \Gamma$ is finite and each $f_{k}$
belongs to $A_{k}$.  Fix $\alpha \in \supp(f)$. Let $n := c(\alpha)$. Since $\alpha \in \supp(f)$, we
have $n \in F$.  That $c$ is continuous and $\Gamma$ is discrete implies $G_n$ is a neighbourhood
of $\alpha$ and $f|_{G_n} = f_n|_{G_n}$ by definition. The function $f_n$ is locally constant, so
there exists a subneighbourhood $U \subseteq G_n$ of $\alpha$ on which $f_n$, and hence $f$, is
constant. Hence $\ag \subseteq \{f \in C_c(\Gg) : f\text{ is locally constant}\}$.

Now suppose that $f \in C_{c}(\Gg)$ is locally constant. Since $c$ is continuous and $\Gp$ is
discrete, each $\Gg_n:=c^{-1}(n)$ is clopen. The collection of open sets $\{\Gg_n\}_{n \in \Gp}$
covers $\supp(f)$ and $\supp(f)$ is compact, so there is a finite subset $F \subseteq \Gp$ such
that $f|_{\Gg_n} \equiv 0$ for $n \not\in F$. Hence, defining $f_n(\alpha) := 1_{\Gg_n}(\alpha)
f(\alpha)$ (multiplication here is pointwise) for all $n \in \Gp$ and $\alpha \in \Gg$, we have $f
= \sum_{n \in F} f_n$. Each $f_n \in A_n$, and  therefore $f \in \ag$. This proves the first
assertion.

For any $U\in\Uu$, the function $1_U$ is locally constant; therefore $\lsp\{1_U : U \in \Uu\}
\subseteq \ag$ by the preceding paragraphs. We must show that $\ag \subseteq \lsp\{1_U : U \in
\Uu\}$. Fix $f \in \ag$. Since $f$ is locally constant and $\Uu$ is a basis, for each $\alpha \in
\supp(f)$, there is a neighbourhood $U_\alpha \in \Uu$ of $\alpha$ such that $f|_{U_{\alpha}}$ is
constant. By Lemma~\ref{lem:clopen support}, we may assume that $U_\alpha \subseteq \supp(f)$.
Since $\supp(f)$ is compact there is a finite subset $F \subseteq \{U_{\alpha}\}_{\alpha \in \supp
f}$ such that $\supp(f) = \bigcup F$. Let $K$ be the disjointification of $F$ discussed in
Remark~\ref{rmk:disjointify}. Since $f$ is constant on each $V \in F$ and each $W \in K$ is a
subset of some $V \in F$, the function $f$ is constant on each $W \in K$. Hence, writing $f(W)$ for
the unique value taken by $f$ on $W \in K$, we have $f = \sum_{W \in K} f(W) 1_{W}$.
\end{proof}

\begin{dfn}\label{dfn:graded-subset}
Let $\Gp$ be a discrete group, $\Gg$ a locally compact, Hausdorff groupoid with totally disconnected
unit space such that $s : \Gg \to \go$ is a local homeomorphism, and $c : \Gg \to \Gp$ a continuous
cocycle. We say that a subset $S$ of $\Gg$ is \emph{graded} if the cocycle $c$ is constant on $S$.
If $S \subseteq c^{-1}(n)$, we say that $S$ is \emph{$n$-graded}. For each $n \in \Gp$ we write
$\BncoG$ for the collection of all $n$-graded compact open bisections of $\Gg$. We write $\BcoG$
for $\bigcup_{n \in \Gp} \BncoG$.
\end{dfn}

\begin{lemma}\label{lem:alg-element-expression}
Let $\Gp$ be a discrete group, $\Gg$ a locally compact, Hausdorff groupoid with totally disconnected
unit space such that $s : \Gg \to \go$ is a local homeomorphism, and $c : \Gg \to \Gp$ a continuous
cocycle. Every $f \in \ag$ can be expressed as $f = \sum_{U \in F} a_U 1_{U}$ where $F$ is a finite
subset of $\BcoG$ whose elements are mutually disjoint and $a : U \mapsto a_U$ is a function from
$F$ to $\CC$.
\end{lemma}
\begin{proof}
Let $f \in \ag$.  By Lemma~\ref{lem:span}, there is a finite set $K_0$ of compact open sets and an
assignment $W \mapsto d_W$ of scalars to the elements of $K_0$ such that $f = \sum_{W \in K_0} d_W
1_W$. Let
\[
    K := \{W \cap G_n : W \in K_0, n \in \Gp, W \cap G_n \not= \emptyset\}.
\]
Since $\Gamma$ is discrete and $c$ is continuous, each $G_n$ is open. Since each $W \in K_0$ is
compact, $K$ is finite. Each $V \in K$ is
 graded; we write $c(V)$ for the unique value taken by $c$ on $V$. For
each $V \in K$, let
\[
b_V = \sum_{W \in K_0, W \cap G_{c(V)} = V} d_W
\]
Then $f = \sum_{V \in K} b_V 1_V$.

Let $F$ be the disjointification of $K$. Each $U \in F$ is graded because $F$ is a refinement of
$K$. For $U \in F$, define
\[
    a_U = \sum_{V \in K,\, U \subseteq V} b_V.
\]
Then $f = \sum_{U \in F} a_U 1_U$ is the desired expression.
\end{proof}

Recall that given a locally compact, Hausdorff groupoid $\Gg$ such that $s : \Gg \to \go$ is a local
homeomorphism, and given $f,g \in A(G) \subseteq C_c(\Gg)$, the functions $f^*$ and $f * g$ are
given by
\begin{align}
    f^*(\gamma) &= \overline{f(\gamma^{-1})} \label{eq:multiplication}\\
    (f*g)(\gamma) &= \sum_{r(\alpha) = r(\gamma)} f(\alpha) g(\alpha^{-1}\gamma).\label{eq:involution}
\end{align}

\begin{prop}\label{prop:A is graded algebra}
Let $\Gp$ be a discrete group, $\Gg$ a locally compact, Hausdorff groupoid with totally
disconnected unit space such that $s : \Gg \to \go$ is a local homeomorphism, and $c : \Gg \to \Gp$
a continuous cocycle. Under the operations \eqref{eq:multiplication}~and~\eqref{eq:involution},
$\ag$ is a $\Gp$-graded $*$-algebra.
\end{prop}

\begin{rmk}
For us, an involution on a $*$-algebra over $\CC$ is always \emph{conjugate} linear.
\end{rmk}

\begin{rmk}  We do not assume that $\Gamma$ is abelian so we will write the group operation
multiplicatively.
\end{rmk}

\begin{proof} That $A(G)$ is a complex algebra follows from \cite[Proposition~4.6]{Steinberg2010}.
We must verify that $A(G)$ is a $*$-algebra and that $A(G)$ is graded. Observe that the $*$-operation is a
conjugate-linear involution on $\ag$ and takes $A_n$ to $A_{n^{-1}}$. Next we will show that the
multiplication defined on $\ag$ is a graded multiplication. If $f \in A_m$ and $g \in A_n$, then if
$(f * g)(\gamma) \not= 0$ we have $f(\alpha) \not= 0$ and $g(\alpha^{-1}\gamma) \not= 0$ for some
$\alpha$ with $r(\alpha) = r(\gamma)$. In particular, $c(\alpha) = m$, and $c(\alpha^{-1}\gamma) =
n$ forcing $c(\gamma) = mn$ (because
$c(\gamma)=c(\alpha\alpha^{-1}\gamma)=c(\alpha)c(\alpha^{-1}\gamma)$). Hence $\supp(f*g)\subseteq
\Gg_{mn}$.
\end{proof}

We finish this section by presenting of $\ag$ as a universal algebra.

\begin{dfn}\label{dfn:representation}
Let $\Gp$ be a discrete group, $\Gg$ a locally compact, Hausdorff groupoid with totally
disconnected unit space such that $s : \Gg \to \go$ is a local homeomorphism, and $c : \Gg \to \Gp$
a continuous cocycle.  Let $B$ be an algebra over $\CC$. A \emph{representation of $\BcoG$} in $B$
is a family $\{t_U : U \in \BcoG\} \subseteq B$ satisfying
\begin{enumerate}\renewcommand{\theenumi}{R\arabic{enumi}}
\item\label{it:zero} $t_\emptyset = 0$;
\item\label{it:multiplicative} $t_Ut_V = t_{UV}$ for all $U,V \in \BcoG$; and
\item\label{it:additive} $t_U + t_V = t_{U \cup V}$ whenever $U$ and $V$ are disjoint elements
    of $\BncoG$ for some $n$.
\end{enumerate}
\end{dfn}

\begin{theorem}\label{thm:UniversalProp}
Let $\Gp$ be a discrete group, $\Gg$ a locally compact, Hausdorff groupoid with totally
disconnected unit space such that $s : \Gg \to \go$ is a local homeomorphism, and $c : \Gg \to \Gp$
a continuous cocycle. Then $\{1_U : U \in \BcoG\} \subseteq \ag$ is a representation of $\BcoG$
which spans $\ag$. Moreover, $\ag$ is universal for representations of $\BcoG$ in the sense that
for every representation $\{t_U : U \in \BcoG\}$ of $\BcoG$ in an algebra $B$, there is a unique
homomorphism $\pi : \ag \to B$ such that $\pi(1_U) = t_U$ for all $U \in \BcoG$.
\end{theorem}
\begin{proof}
The collection $\{1_U : U \in \BcoG\}$ certainly satisfies (\ref{it:zero})~and~(\ref{it:additive}),
and it satisfies~(\ref{it:multiplicative}) by \cite[Proposition~4.5 (3)]{Steinberg2010}. That this family spans
$\ag$ follows from Lemma~\ref{lem:alg-element-expression}.

Let $B$ be a complex algebra and let $\{t_U : U \in \BcoG\}$ be a representation of $\BcoG$ in $B$.
We must show that there is a homomorphism $\pi : \ag \to B$ satisfying $\pi(1_U) = t_U$ for all $U
\in \BcoG$; uniqueness follows from the previous paragraph. We begin by showing that
\begin{equation}\label{eq:old additive}
\parbox{0.8\textwidth}{
    $\vphantom{\Big|}\sum_{U \in F} t_U = t_{\bigcup F}$ for $n \in \Gp$ and finite $F \subseteq \BncoG$
consisting of mutually disjoint bisections such that $\bigcup F \in \BncoG$.
}
\end{equation}
Let $F\subseteq \Bco{n}{G}$ be a finite collection of mutually disjoint bisections
such that $\bigcup F$ is a bisection. We claim that $r(U) \cap r(V) = \emptyset$ for distinct $U,V
\in F$. To see this, fix $x \in r(U)$. There exists $\alpha \in U$ such that $r(\alpha) = x$, and
this $\alpha$ is the unique element of $\bigcup F$ whose range is $x$ because $\bigcup F$ is a
bisection. Since $U \cap V = \emptyset$, we have $\alpha \not\in V$ and hence $x \not\in r(V)$. So
the sets $r(U)$ where $U \in F$ are mutually disjoint as claimed. Thus each $U \in F$ satisfies $U
= r(U) (\bigcup F)$. A standard induction extends~(\ref{it:additive}) to finite collections of
mutually disjoint compact open subsets of $\go$. Combining this with~(\ref{it:multiplicative}), we
obtain
\[
t_{\bigcup F}
     = t_{r(\bigcup F)} t_{\bigcup F}
     = \sum_{U \in F} t_{r(U)} t_{\bigcup F}
     = \sum_{U \in F} t_{r(U) \big(\bigcup F\big)}
     = \sum_{U \in F} t_U.
\]

We show next that the formula $\sum_{U \in F} a_U 1_U \mapsto \sum_{U \in F} a_U t_U$ is
well-defined on linear combinations of indicator functions where $F \subseteq \BcoG$ is a finite
collection of mutually disjoint bisections. It will follow from Lemma
\ref{lem:alg-element-expression} that there is a unique linear map $\pi : A(G) \to B$
such that $\pi(1_U) = t_U$ for each $U \in \BcoG$. Fix $f \in \ag$ and suppose that
\[\textstyle
    \sum_{U \in F} a_U 1_U = f = \sum_{V \in H} b_V 1_V
\]
where each of $F$ and $H$ is a finite set of mutually disjoint elements of $\BcoG$. We must show
that
\[\textstyle
    \sum_{U \in F} a_U t_U = \sum_{V \in H} b_V t_V.
\]
Since the $\Gg_n$ are mutually disjoint, for each $n \in \Gp$ we
have
\[
\sum_{U \in F \cap \BncoG} a_U 1_U = f|_{G_n} = \sum_{V \in H \cap \BncoG} b_V 1_V,
\]
so we may assume that $F, G \subseteq \BncoG$ for some $n \in \Gp$.

Let $K = \{U \cap V : U \in F, V \in H, U \cap V \not= \emptyset\}$. Then each $W \in K$ belongs to
$\BncoG$. Moreover, for $U \in F$ we have $U = \bigsqcup \{W \in K: W \subseteq U\}$.
Hence~\eqref{eq:old additive} gives $t_U = \sum_{W \in K, W \subseteq U} t_W$ for each $U \in F$; a
similar decomposition holds for $t_V$ for each $V \in H$. Therefore
\[
    \sum_{U \in F} a_U t_U
        = \sum_{U \in F}\  \sum_{W \in K,\,W \subseteq U} a_U t_W
        = \sum_{W \in K} \Big(\sum_{U \in F,\, W \subseteq U} a_U\Big) t_W,
\]
and similarly
\[
    \sum_{V \in F}b_V t_V
        = \sum_{W \in K} \Big(\sum_{V \in F,\, W \subseteq V} b_V\Big) t_W.
\]
Fix $W \in K$. It suffices now to show that $\sum_{U \in F,\, W \subseteq U} a_U = \sum_{V \in F,\,
W \subseteq V} b_V$. By definition of $K$, the set $W$ is nonempty, so let $\alpha \in W$. Then
for $U \in F$, we have $\alpha \in U \implies W \cap U \not= \emptyset \implies W \subseteq
U$. Since $\alpha \in W$, this implies that $\alpha \in U \iff W \subseteq U$. Hence
\[
f(\alpha) = \sum_{U \in F} a_U 1_U(\alpha)
    = \sum_{U \in F,\, \alpha \in U} a_U
    = \sum_{U \in F,\, W \subseteq U} a_U.
\]
a similar calculation shows that $\sum_{V \in F,\, W \subseteq V} b_V = f(\alpha)$ as
well. So there is a linear map $\pi : A(G) \to B$ such that
$\pi(1_U) = t_U$ for all $U \in \BcoG$.

We must check that $\pi$ is a homomorphism.  To see that $\pi$ is multiplicative, fix $f, g \in
\ag$. Express $f = \sum_{U \in F} a_U 1_U$ and $g=\sum_{V \in H} b_V 1_V$ where $F$ and $H$ are finite
subsets of $\BcoG$, and calculate:
\[
\pi(fg)
    = \pi\Big(\Big(\sum_{U \in F} a_U 1_U\Big)\Big(\sum_{V \in H} b_V 1_V\Big)\Big)\\
    = \pi\Big(\sum_{U \in F}\sum_{V \in H} a_U b_V 1_{U} 1_{V}\Big).
\]
Since \cite[Proposition~4.5 (3)]{Steinberg2010} gives $1_{U} 1_{V} = 1_{UV}$ for all $U,V$, we then have
\[
\pi(fg)
    = \pi\Big(\sum_{U \in F}\sum_{V \in H} a_U b_V 1_{UV}\Big)
    = \sum_{U \in F}\sum_{V \in H} a_U b_V t_{UV}.
\]
Each $t_{UV} = t_{U} t_{V}$ by~(\ref{it:multiplicative}), so
\[
\pi(fg)
    = \sum_{U \in F}\sum_{V \in H} a_U b_V t_{U}t_{V}
    = \Big(\sum_{U \in F} a_U t_U\Big)\Big(\sum_{V \in H} b_V t_V\Big)
    = \pi(f)\pi(g)
\]
as required.
\end{proof}

\section{\texorpdfstring{$\ag$}{A(G)} is dense in
\texorpdfstring{$\cs(\Gg)$}{C*(G)}}\label{sec:dense}

Since our aim is to produce algebras associated to totally disconnected, locally compact, Hausdorff
groupoids whose relationship to the groupoid {$\cs$-}algebra is analogous to that of Leavitt path
algebras to graph {$\cs$-}algebras, we show in this section that the subalgebra $\ag$ of $C_c(\Gg)$
is dense in the full (and hence also the reduced) $C^*$-algebra of $\Gg$. We could prove this as in
\cite[Proposition~4.1]{KumjianPaskEtAl:JFA97} or \cite[Proposition~6.7]{Steinberg2010}
by using the Stone-Weierstrass theorem, but a direct argument takes about the same amount of effort.

We first prove a technical lemma.

\begin{lemma}\label{lem:cover}
Let $\Gg$ be a locally compact, Hausdorff groupoid with
totally disconnected unit space such that $s : \Gg \to \go$ is a local
homeomorphism. Fix a compact open
bisection $U$ of $\Gg$ and suppose that $f \in C_c(\Gg)$ is supported on $U$. Fix $\varepsilon
> 0$. There exists a finite set $\mathcal{V}$ of nonempty compact open bisections of $\Gg$ such
that $U = \bigsqcup \mathcal{V}$ and such that for each $V \in \mathcal{V}$, we have
$|f(\alpha) - f(\beta)| \le \varepsilon$ for all $\alpha, \beta \in V$.
\end{lemma}
\begin{proof}
For each $\gamma \in U$ let $U_\gamma$ be a compact open neighbourhood of $\gamma$ such that
$U_\gamma \subseteq U$ and $|f(\alpha) - f(\gamma)| < \varepsilon/2$ for all $\alpha \in U_\gamma$.
Since $U$ is compact, there is a finite subset $F$ of $U$ such that $\{U_\gamma : \gamma \in F\}$
covers $U$. Let $\mathcal{V}$ be the disjointification of the $U_\gamma$ as in
Remark~\ref{rmk:disjointify}. Fix $V\in\mathcal{V}$. Then there exists $\gamma \in F$ such that $V
\subseteq U_\gamma$, and then for $\alpha, \beta \in V$, we have $|f(\alpha) - f(\beta)| \le
|f(\alpha) - f(\gamma)| + |f(\gamma) - f(\beta)| < \varepsilon$.
\end{proof}

To state the next proposition, we recall from \cite{Ren} that for a locally compact Hausdorff
groupoid $G$ in which $s : \Gg \to \go$ is a local homeomorphism,
the $I$-norm on $C_c(\Gg)$ is defined as follows.  For $f \in
C_c(\Gg)$, let
\[
\|f\|_{I,r}:= \sup_{u \in \go} \Big\{ \sum_{r(\alpha)=u}|f(\alpha)|\Big\}\quad
\text{and} \quad \|f\|_{I,s}:= \sup_{u \in \go}\Big\{\sum_{s(\alpha)=u}|f(\alpha)|\Big\}.
\]
Then the $I$-norm of $f$ is $\|f\|_I:=\max\{\|f\|_{I,r}, \|f\|_{I,s}\}$.  The $I$-norm
dominates each of the universal norm, the reduced norm, and the uniform norm on $C_c(\Gg)$. (See
\cite{Ren} for further details.)

\begin{prop}\label{prop:AG dense in CcG}
Let $\Gp$ be a discrete group, $\Gg$ a locally compact, Hausdorff groupoid with totally disconnected
unit space such that $s : \Gg \to \go$ is a local homeomorphism, and $c : \Gg \to \Gp$ a continuous
cocycle. With notation as above, $\ag$ is dense in $C_c(\Gg)$ under each of the reduced norm, the
universal norm, the $I$-norm, and the uniform norm.
\end{prop}

\begin{proof}
Since the $I$-norm dominates the other  three norms, it suffices to prove the result for the
$I$-norm. Fix $f \in C_c(\Gg)$ and $\varepsilon > 0$. Since $f$ has compact support, $\supp (f)$
can be written as a finite union of elements of $\BcoG$. So we can write $f =\sum_{i=1}^n f_i$
where each $f_i$ is supported on an element of $\BcoG$.  For each $i$, apply Lemma~\ref{lem:cover}
to $\supp(f_i)$ to obtain a cover $\Uu_i$ of the support of $f_i$ by disjoint compact open
bisections such that for $U \in \Uu_i$, we have $|f_i(\alpha) - f_i(\beta)| \leq \varepsilon/n$ for
all $\alpha, \beta \in U_i$. For each $i \le n$ and each $U \in \Uu_i$, fix $z_{i,U} \in f(U)$, so
$|f_i(\alpha) - z_{i,U}| \leq \varepsilon/n$ for all $\alpha \in U$.  Then let $g_i := \sum_{U
\in \Uu_i} z_{i,U} 1_U$ for all $i\leq n$ and define $g := \sum_{i=1}^n g_i \in \ag$. We have
\[
   \|f - g\|_I \le \sum_{i=1}^n \|f_i - g_i\|_I.
\]
Fix $i \le n$. It suffices to show that $\|f_i - g_i\|_I \le \varepsilon/n$. Fix $u \in \go$. Since
$f_i$ is supported on a bisection, there is at most one $\alpha \in s^{-1}(u) \cap \supp(f_i)$. If
there is no such $\alpha$, then $\sum_{s(\alpha)=u} |(f_i - g_i)(\alpha)| = 0$ and we are done. So
suppose that $\alpha \in s^{-1}(u) \cap \supp(f_i)$. Then there is a unique $U_0 \in \Uu_i$ such
that $\alpha \in U_0$.  Therefore $\sum_{s(\alpha)=u} |(f_i - g_i)(\alpha)| = |f_i(\alpha) - z_{i,U_0}|
\leq \varepsilon/n$. Since $u \in \go$ was arbitrary, we conclude that $\|f_i - g_i\|_{I, s} \le
\varepsilon/n$. A symmetric argument gives $\|(f_i - g_i)(\alpha)\|_{I, r} \leq \varepsilon/n$. Hence
$\|f_i - g_i\|_I  \leq \varepsilon / n$ as required.
\end{proof}

\begin{prop}\label{prop:agL isomorphic to KP algebra}
Suppose that $\Lambda$ is a row-finite, $k$-graph with no sources and that $\Gg_\Lambda$ is the
corresponding $k$-graph groupoid. Then $\agL$ as constructed above is isomorphic to the
Kumjian-Pask algebra $\operatorname{KP}(\Lambda,\CC)$.
\end{prop}
\begin{proof}
By \cite[Corollary 3.5]{KumjianPask:NYJM00},  $t_\lambda :=
1_{Z(\lambda, s(\lambda))}$ determines a Cuntz-Krieger $\Lambda$-family in $C^*(\Gg)$. In
particular, there is a Kumjian-Pask family (\cite[Definition 3.1]{ArandaClarkEtAl:xx11}) for
$\Lambda$ determined by $t_\lambda = 1_{Z(\lambda, s(\lambda))}$ and $t_{\lambda^*} =
1_{Z(s(\lambda),\lambda)}$ for all $\lambda \in \Lambda$.  It follows from the universal property
of $\operatorname{KP}(\Lambda,\CC)$ that there is a homomorphism $\phi :
\operatorname{KP}(\Lambda,\CC) \to \agL$ which carries each $s_\lambda$ to $t_\lambda$ and each
$s_{\lambda^*}$ to $t_{\lambda^*}$.

By \cite[Theorem 3.4]{ArandaClarkEtAl:xx11} the algebra $\operatorname{KP}(\Lambda, \CC)$ is
spanned by the elements $t_\mu t_{\nu^*}$ where $\mu,\nu \in \Lambda$ with $s(\mu) = s(\nu)$, and
the $\ZZ$-grading of $\operatorname{KP}(\Lambda,\CC)$  carries each $s_\mu s_{\nu^*}$ to
$d(\mu)-d(\nu)$.  So to see that $\phi $ is graded, it suffices to show that it preserves the
grading of each $s_\mu s_{\nu^*}$, which it does since
\[
    \phi(s_\mu s_{\nu^*}) = 1_{Z(\mu, \nu)} = 1_{\{(\mu x, d(\mu)-d(\nu), \nu x) : x \in \Lambda^\infty, r(x)=s(\mu)}\}
\in A_{d(\mu)-d(\nu)}.
\]
Since each $Z(v)$ is nonempty, $\phi(p_v) \neq 0$ for each $v \in E^0$.  Thus the graded uniqueness
theorem for Kumjian-Pask algebras \cite[Theorem 4.1]{ArandaClarkEtAl:xx11}  implies that $\phi$ is
injective.

It remains to show that $\phi$ is surjective. By Lemma~\ref{lem:span}, $\agL$ is spanned by the
functions $1_U$ where $U$ ranges over  all compact open bisections in $\Gg_\Lambda$.  Let $U$ be a
compact open bisection. Since the grading is continuous and $U$ is compact, we can write $1_U$ as
the finite sum $\sum_{U \cap G_n \not= \emptyset} 1_{U \cap G_n}$ where each $U \cap G_n$ is a
graded compact open bisection. So fix $n \in \NN^k$ and a compact open $n$-graded bisection $V$. It
suffices to show that $1_V \in \lsp\{1_{Z(\mu,\nu)} : s(\mu) = s(\nu)\}$. Because $V$ is compact
and the sets $Z(\mu,\nu)$ form  a basis for the topology on $G_{\Lambda}$
\cite[Proposition~2.8]{KumjianPask:NYJM00}, we can write $V = \bigcup_{(\mu,\nu) \in F} Z(\mu,\nu)$
for some finite set $F \subseteq \{(\mu, \nu) \in \Lambda \times \Lambda: s(\mu) = s(\nu)\}$. Since
$V$ is $n$-graded, we have $d(\mu) - d(\nu) = n$ for all $(\mu,\nu) \in F$. Let $p :=
\bigvee_{(\mu,\nu) \in F} d(\mu)$. Then for each $(\mu,\nu) \in F$ we have $Z(\mu,\nu) = \bigcup
\{Z(\mu\alpha, \nu\alpha) : \alpha \in s(\mu)\Lambda^{p - d(\mu)}\}$. Let ${H} :=
\{(\mu\alpha,\nu\alpha) : (\mu,\nu) \in F, \alpha \in s(\mu)\Lambda^{p - d(\mu)}\}$. Then
$Z(\eta,\zeta) \cap Z(\eta',\zeta') = \emptyset$ for distinct $(\eta,\zeta), (\eta',\zeta') \in
{H}$, so $V = \bigsqcup_{(\eta,\zeta) \in {H}} Z(\eta,\zeta)$. Hence $1_U = \sum_{(\eta, \zeta) \in
{H}} 1_{Z(\eta, \zeta)}$, and it follows that $\phi$ is surjective.
\end{proof}

\begin{rmk}
When $k = 1$ in the preceding proposition, $\Lambda$ is the path category of the directed graph $E
= (\Lambda^0, \Lambda^1, r, s)$ and, in this case, the proposition specialises to the statement that
$\ag$ is isomorphic to the Leavitt path algebra of \cite{AbrPino}.
\end{rmk}

\section{The uniqueness theorems}\label{sec:uniqueness}

Interestingly, in the situation of groupoids, the graded uniqueness theorem is a corollary of the
natural generalisation of the Cuntz-Krieger uniqueness theorem. This in turn is essentially
Renault's structure theorem for the reduced $C^*$-algebra of a groupoid in which the units with
trivial isotropy are dense in the unit space.

\begin{rmk}
The condition that the units with trivial isotropy are dense in the unit space occurs frequently in
the groupoid literature.  It has been variously referred to as ``topologically free,''
``essentially principal,'' and ``essentially free.'' It seems that
``topologically free" is becoming the standard term, but since we use the hypothesis
only in a few places, we avoid any confusion by stating it in full each time.
\end{rmk}

In the groupoid literature, given a unit $u$, it is standard to denote the isotropy subgroup
$\{\alpha \in G : r(\alpha) = s(\alpha) = u\}$ by either $G(u)$ or $G^u_u$. Here we have chosen the
more suggestive notation $uGu$, which is in keeping with the notation established
in~\eqref{eq:setproduct}. Likewise, we write $\Gg u$ for $s^{-1}(u)$.

\begin{theorem}\label{thm:CKUT}
Let  $\Gg$ be a locally compact, Hausdorff groupoid with totally disconnected unit space such that
$s : \Gg \to \go$ is a local homeomorphism.  Suppose that $\{u \in \go : u\Gg u = \{u\}\}$ is dense
in $\go$. Let $\pi : \ag \to B$ be a homomorphism. Suppose that $\ker(\pi) \not= \{0\}$. Then there
is a compact open subset $K \subseteq \go$ such that $\pi(1_K) = 0$.
\end{theorem}

To prove our theorem, we  need a technical lemma.

\begin{lemma}\label{lem:bisection neighborhood}
Let $G$ be a locally compact, Hausdorff groupoid such that $s: G \to \go$ is a local homeomorphism.
Fix $\alpha \in \Gg$ and a precompact neighbourhood $V$ of $\alpha$. Suppose that $r(\alpha) \Gg
s(\alpha) = \{\alpha\}$. Then there exist neighbourhoods $X$ of $r(\alpha)$ and $Y$ of $s(\alpha)$
such that $X V Y$ is a precompact open bisection.
\end{lemma}
\begin{proof}
Suppose, to the contrary, that for every neighbourhood $X$ of $r(\alpha)$ and every neighbourhood
$Y$ of $s(\alpha)$, $XVY$ fails to be a bisection.  Let $U$ be an open bisection containing
$\alpha$. Fix a fundamental sequence of neighbourhoods
$(Y_i)^\infty_{i=1}$ of $s(\alpha)$, and for each $i$, let $X_i := r(U Y_i)$, so that
$(X_i )_{i=1}^\infty$ forms a fundamental sequence of
neighbourhoods of $r(\alpha)$. Since each $X_i V Y_i$ fails to be a bisection, for each $i$ there
exist $\beta_i, \gamma_i \in X_i V Y_i$ with $\beta_i \neq \gamma_i$ such that either $s(\beta_i) =
s(\gamma_i)$ or $r(\beta_i) = r(\gamma_i)$ for all $i$. The sequence $\big((\beta_i,
\gamma_i)\big)^\infty_{i=1}$ belongs to the precompact set $V \times V$, so by passing to a
subsequence and relabelling we may assume that $\beta_i \to \beta$ and $\gamma_i \to \gamma$.
Since the $X_i$ and $Y_i$ are fundamental sequences of neighbourhoods,
it follows that $r(\beta_i), r(\gamma_i) \to r(\alpha)$ and $s(\beta_i), s(\gamma_i) \to
s(\alpha)$. Since $r,s : \Gg \to \go$ are continuous and $\go$ is Hausdorff, $r(\beta) = r(\alpha)
= r(\gamma)$ and $s(\beta) = s(\alpha) = s(\gamma)$. By hypothesis, $s(\alpha) \Gg r(\alpha) =
\{\alpha\}$, so we have $\beta = \gamma = \alpha$. Since $U$ is a neighbourhood of $\alpha$, we
then have $\beta_i, \gamma_i \in U$ for large $i$.  Fix $i$ such that $\beta_i, \gamma_i \in U$.
Then $\beta_i \not = \gamma_i$ but either $r(\beta_i) = r(\gamma_i)$ or $s(\beta_i) = s(\gamma_i)$,
contradicting that $U$ is a bisection.
\end{proof}

\begin{proof}[Proof of Theorem~\ref{thm:CKUT}]
Fix $f \in \ker(\pi) \setminus \{0\}$. Since $s$ is a local homeomorphism, it is an open map, so
$s(\supp(f)) \subseteq \go$ is open by Lemma~\ref{lem:clopen support}.  Because $\{u \in \go : u\Gg u
= \{u\}\}$ is dense in $\go$, there exists $u \in s(\supp(f))$ such that $u\Gg u = \{u\}$. Fix
$\alpha \in \supp(f)$ with $s(\alpha) = u$. Then $r(\alpha)\Gg s(\alpha) = \alpha (\alpha^{-1} \Gg
u) \subseteq \alpha (u\Gg u) = \{\alpha\}$.

By Lemma~\ref{lem:bisection neighborhood}, there exist compact open neighbourhoods $X$ of
$r(\alpha)$ and $Y$ of $s(\alpha)$ such that $X\supp(f)Y$ is a bisection containing $\alpha$.
Because $r$ and $s$ are continuous, $X \supp(f) Y=  r^{-1}(X) \cap \supp(f) \cap s^{-1}(Y)$ is
compact and by Lemma \ref{lem:clopen support} it is also open. Since $f$ is locally constant, there
exist subneighbourhoods $X_0 \subseteq X$ of $r(\alpha)$ and $Y_0 \subseteq Y$ of $s(\alpha)$ such
that $X_0 \supp(f) Y_0$ is a compact open bisection and $f(\beta) = f(\alpha)$ for all $\beta \in
X_0 \supp(f) Y_0$.

We have $1_{X_0}, 1_{Y_0} \in \ag$. By Lemma~\ref{lem:span}, $f$ may be written as a linear
combination of characteristic functions of compact open bisections.
\cite[Proposition~4.5 (3)]{Steinberg2010} together with bilinearity of multiplication implies that for $\beta \in \Gg$,
\[
    (1_{X_0} * f * 1_{Y_0})(\beta)
        = 1_{X_0}(r(\beta)) f(\beta) 1_{X_0}(s(\beta))
        = 1_{X_0 \supp(f) Y_0}(\beta) f(\beta)
        = 1_{X_0 \supp(f) Y_0}(\beta) f(\alpha).
\]
Thus $f_0 := 1_{X_0} * f * 1_{Y_0} = f(\alpha) 1_{X_0 \supp(f) Y_0}$. Since $\pi(f)=0$, we have
$\pi(f_0) = 0$. We have $(X_0 \supp(f) Y_0)^{-1}(X_0 \supp(f) Y_0) = Y_0$ because $X_0 \supp(f)
Y_0$ is a bisection. \cite[Proposition~4.5 (3)]{Steinberg2010} implies that
\[
f^*_0 * f_0 = |f(\alpha)|^2 1_{(X_0 \supp(f) Y_0)^{-1}(X_0 \supp(f) Y_0)} = |f(\alpha)|^2 1_{Y_0}.
\]
Hence $K := Y_0$ satisfies $\pi(1_K) =
\frac{1}{|f(\alpha)|^2} \pi(f^*_0
* f_0) = 0$ as required.
\end{proof}

Our graded uniqueness theorem now follows from a bootstrapping argument.

\begin{theorem}\label{thm:kernel of graded homomorphism}
Let $\Gp$ be a discrete group, $\Gg$ a locally compact, Hausdorff groupoid with totally disconnected
unit space such that $s : \Gg \to \go$ is a local homeomorphism, and $c : \Gg \to \Gp$ a continuous
cocycle. Suppose that $\{u \in \go : u\Gg_e u = \{u\}\}$ is dense in $\go$. Let $\pi : \ag \to B$
be a graded homomorphism, and suppose that $\ker(\pi) \not= \{0\}$. Then there is a compact open
subset $K \subseteq \go$ such that $\pi(1_K) = 0$.
\end{theorem}
\begin{proof}
We first claim that there exists nonzero $f \in A_e$ such that $\pi(f) = 0$. To see this, observe
that since $\ker (\pi) \neq 0$, there exists $g \in \ker(\pi)\setminus\{0\}$ such that $\pi(g) =
0$. Since $g$ is an element of the graded algebra $A(G)$, $g$ can be expressed as a finite sum of
graded components $g = \sum_{h \in F}  g_h$ where $F \subseteq \Gamma$ and each $g_h \in A_h$. Now
$\pi(g) = \sum_{h \in F} \pi(g_h) = 0$, and each $\pi(g_h) \in B_h$ because $\pi$ is a graded
homomorphism. Because the graded subspaces of $B$ are linearly independent, it follows that each
$\pi(g_h) = 0$. Since $g \not= 0$, there exists $k \in F$ such that $g_{k} \not= 0$.  By
Lemma~\ref{lem:alg-element-expression}, we can write $g_{k}$ as $\sum_{V \in K} a_V 1_{V}$ where
$K$ is a finite set of mutually disjoint elements of $\Bco{k}{\Gg}$.  Note that $g_k^* = \sum_{V \in K}
\overline{a_V} 1_{V^{-1}}$; define $f := g_k^* * g_k$. We claim that $f \in A_e \setminus
\{0\}$ and $\pi(f) = 0$. To see this, first notice that
\[
f = \Big(\sum_{V \in K} \overline{a_V} 1_{V^{-1}}\Big) * \Big(\sum_{W \in K} a_W 1_{W}\Big)
    =  \sum_{V,W \in K} \overline{a_V}a_W 1_{V^{-1}} * 1_{W}
    =  \sum_{V,W \in K} \overline{a_V}a_W 1_{V^{-1}W}
\]
by~\cite[Proposition~4.5 (3)]{Steinberg2010}.  Now, because each $V \in K$ is a subset of $G_{k}$, each $V^{-1}W
\subseteq G_{k^{-1} k} = G_e$, and thus $f \in A_e$ as claimed. We have $\pi(f) = 0$ because
$\pi(g_{k}) = 0$.

To show that $f$ is nonzero, fix $\alpha \in G_k$ such that $g(\alpha) \not= 0$. Since the elements
of $K$ are mutually disjoint, there is a unique $V_\alpha \in K$ such that $\alpha \in V_\alpha$,
and then $a_{V_\alpha} = g(\alpha) \not= 0$. Since $s$ is a local homeomorphism, $\Gg s(\alpha)$ is
a discrete space. Write $C_c(G s(\alpha))$ for the space of finitely supported functions from
$G s(\alpha)$ to $\CC$ and for each $\beta \in \Gg s(\alpha)$ let $\delta_\beta$ denote the
point-mass at $\beta$ so that $C_c(\Gg s(\alpha)) = \lsp\{\delta_\beta : \beta \in \Gg
s(\alpha)\}$. For $f \in C_c(G)$, let $\rho(f)$ be the linear map on $C_c(G s(\alpha))$ determined
by
\[
\rho(f) \delta_{\beta} = \sum_{s(\alpha) = r(\beta)} f(\alpha) \delta_{\alpha\beta}.
\]
Let $(\cdot\,|\,\cdot)$ be the standard inner product on $C_c(G s(\alpha))$, that is $(f|g) =
\sum_\beta \overline{f(\beta)} g(\beta)$. Since the elements of $K$ are mutually disjoint,
$\big(\rho(1_V) \delta_{s(\alpha)} \mid \rho(1_W) \delta_{s(\alpha)}\big) = 0$ for distinct $V,W
\in K$. A calculation shows that for $V \in K$ and $\beta, \gamma \in G s(\alpha)$, we have
$(\delta_\beta | \rho(1_{V^{-1}}) \delta_\gamma) = (\rho(1_V) \delta_\beta | \delta_\gamma)$. Hence
\begin{align*}
\big(\rho(f) \delta_{s(\alpha)} \mid \delta_{s(\alpha)}\big)
    &= \big(\rho(g_k) \delta_{s(\alpha)} \mid \rho(g_k) \delta_{s(\alpha)}\big) \\
    &= \sum_{V,W \in K} \overline{a_V} a_W
                                            \big(\rho(1_W) \delta_{s(\alpha)} \mid \rho(1_V) \delta_{s(\alpha)}\big)
    = \sum_{\substack{V \in K,\\ s(\alpha) \in s(V)}} |a_V|^2
    \ge |a_{V_\alpha}|^2.
\end{align*}
Hence $\rho(f) \not= 0$ which forces $f \not= 0$.

By hypothesis  $\{u\in \go: u\Gg_e u=\{u\}\}$ is dense in $\go$. By definition, $A_e$ is equal
to the space of locally constant, continuous, compactly supported functions on $\Gg_e$, so we may
apply Theorem~\ref{thm:CKUT} to see that $\pi|_{A_{e}} : A_e \to B$ annihilates $1_K$ for some
compact open $K \subseteq \go_e = \go$.
\end{proof}

\begin{cor}\label{cor:gut} Let $\Gp$ be a discrete group, $\Gg$ a locally
compact Hausdorff groupoid with totally disconnected unit space such that $s : \Gg \to \go$ is a
local homeomorphism, and $c : \Gg \to \Gp$ a continuous cocycle. Suppose that $\{u \in \go : u\Gg_e
u = \{u\}\}$ is dense in $\go$. Let $B$ be a $\Gp$-graded complex algebra and let $\{t_U : U \in
\BcoG\}$ be a representation of $\BcoG$ in $B$. Suppose that $t_U \in B_n$ whenever $U \in \BncoG$
and that $t_K \not= 0$ for each compact open $K \subseteq \go$. Then the homomorphism $\pi : \ag
\to B$ obtained from Theorem~\ref{thm:UniversalProp} is injective.
\end{cor}
\begin{proof}
Since each $\ag_n$ is spanned by $\{1_{U} : U \in \BncoG\}$, the homomorphism $\pi$ is graded.
Since $\pi(1_K) = t_K \not= 0$ for all compact open $K \subseteq \go$, it follows from
Theorem~\ref{thm:kernel of graded homomorphism} that $\ker(\pi) = \{0\}$.
\end{proof}

\begin{rmk}\label{rmk:ckut rephrase} Suppose that $\Gg$ is a locally compact
Hausdorff groupoid with totally disconnected locally compact unit space such that $s : \Gg \to \go$
is a local homeomorphism and such that $\{u \in \go : u\Gg u = \{u\}\}$ is dense in $\go$. We may
apply Corollary~\ref{cor:gut} with $c$ the trivial cocycle to prove that $\ag$ is the unique
algebra generated by nonzero elements $\{t_U : U \text{ is a compact open bisection of
}\Gg\}$ satisfying
\begin{enumerate}
\item $t_\emptyset = 0$;
\item $t_Ut_V = t_{UV}$ for all compact open bisections $U,V$; and
\item $t_U + t_V = t_{U \cup V}$ whenever $U$ and $V$ are disjoint
    compact open bisections.
\end{enumerate}
\end{rmk}

\begin{rmk}
In the proof of Theorem~\ref{thm:kernel of graded homomorphism}, to see that the function $g_k^* *
g_k$ was nonzero, we really just checked that its image under Renault's left-regular
representation of $\Gg$ associated to the unit $s(\alpha)$ is nonzero.  However, since we are not
working in a $C^*$-completion, we can do everything at the level of linear algebra rather than on
Hilbert space. We could instead have appealed to the $C^*$-identity by regarding $A(\Gg)$ as a
subalgebra of $C_r(\Gg)$, but chose a more elementary argument: our argument is
essentially that used by Renault to show that the reduced norm is positive
definite on $C_c(\Gg)$.
\end{rmk}

\begin{rmk} Recall from \cite{FMY} that if $\Lambda$ is a finitely aligned $k$-graph,
then the $k$-graph groupoid $\Gg_\Lambda$ is totally disconnected and locally compact, and carries
a $\ZZ^k$-grading such that ${\{u \in \go : u\Gg_e u = \{u\}\}}$ is dense in $\go$. So our graded
uniqueness theorem applies to $A(\Gg_\Lambda)$ for any finitely aligned $k$-graph. Likewise,
Remark~\ref{rmk:ckut rephrase} suggests a Cuntz-Krieger uniqueness theorem for $A(\Gg_\Lambda)$.
But in practise the relations described in Definition~\ref{dfn:representation} and
Remark~\ref{rmk:ckut rephrase} are much harder to verify than those of
\cite[Definition~3.1]{ArandaClarkEtAl:xx11}.
\end{rmk}

\begin{example}
To see that the class of Leavitt groupoid algebras is broader than that of Kumjian-Pask algebras
\cite{ArandaClarkEtAl:xx11}, we describe a class of groupoids that satisfy our hypothesis but do
not necessarily arise from $k$-graphs.  Let $T:X\to X$ be a surjective local homeomorphism of a
totally disconnected, compact, Hausdorff space $X$.  Define $T^0 := \id$ and for $k \ge 2$ let $T^k
:= T \circ \cdots \circ T$ be the $k$-fold self-composite of $T$. Let $\Gg$ be the Deaconu-Renault
groupoid defined in \cite[Section~3]{IonescuMuhly}. So
\[
    \Gg = \{(x,n,y)\in X\times\mathbb{Z}\times X\;:\; T^{k}(x)=T^{l}(y),n=k-l\}.
\]
Let $\go$ be the subset $\{(x,0,x) : x \in X\}$, which we identify with $X$ in the obvious way. The
range and source maps are given by $r(x,n,y) = x$ and $s(x,n,y) = y$. Hence triples
$(x_{1},n_{1},y_{1})$ and $(x_{2},n_{2},y_{2})$ are composable if and only if $x_{2}=y_{1}$, in
which case $(x_{1},n_{1},y_{1})(x_{2},n_{2},y_{2}) := (x_{1},n_{1}+n_{2},y_{2})$. The inverse of
$(x,n,y)$ is $(y,-n,x)$. For open subsets $U,V \subseteq X$ and $k,l \ge 0$ such that $T^{k}|_U$
and $T^{l}|_V$ are homeomorphisms and $T^{k}(U)=T^{l}(V)$, define
\[
    Z(U,V,k,l):=\{(x,k-l,y)\in \Gg:x\in U,y\in V\}.
\]
Then
\begin{align*}
    \{Z(U,V,k,l) : U,{}&V \subseteq X\text{ are compact open}, k,l \ge 0,\\
        &T^{k}|_U\text{ and }T^{l}|_V\text{ are homeomorphisms and } T^{k}(U)=T^{l}(V)\}
\end{align*}
is a basis of compact open sets for a topology on $\Gg$ under which it becomes a locally compact,
Hausdorff groupoid with totally disconnected unit space $X$. Fix $(x, n, y) \in \Gg$ and $k,l$ such
that $k - l = n$ and $T^k(x) = T^l(y)$. The source map on $\Gg$ restricts to a homeomorphism on
each basic open set $Z(U,V,k,l)$ so is a local homeomorphism. Moreover, the map $c : \Gg \to
\mathbb{Z}$ defined by $c((x,n,y)) = n$ is a cocycle and is continuous because each basic open set
belongs to some $c^{-1}(n)$. Hence $(\Gg, c)$ satisfies our hypotheses, and $\ag$ is a sensible
candidate for the Leavitt algebra of $(X, T)$.
\end{example}

\end{document}